\documentclass[11pt,A4]{article}
\usepackage{amsmath,amsthm,amsfonts,amssymb,amscd, amsxtra,mathrsfs, color}
\usepackage[english]{babel}
\usepackage{setspace}

\usepackage[margin=1.0in]{geometry}
\newtheorem{theorem}{Theorem}
\newtheorem{lemma}[theorem]{Lemma}

\newtheorem{proposition}[theorem]{Proposition}
\newtheorem{remark}[theorem]{Remark}
\newtheorem{definition}[theorem]{Definition}

\newtheorem{algorithm}{Algorithm}

\begin{document}

\title{A version of the inexact Levenberg-Marquardt method for constrained nonsmooth equations}

\author{Fabiana R. de Oliveira \thanks{Instituto de Matem\'atica e Estat\'istica, Universidade Federal de Goi\'as,  CEP 74001-970 - Goi\^ania, GO, Brazil, E-mails: {\tt  rodriguesfabiana@ufg.br.com}, {\tt  fabriciarodrigues@ufg.br}. }
\and
Fabr\'icia R. Oliveira \footnotemark[1]
}

\maketitle
\begin{abstract}
We herein propose a variant of the projected inexact Levenberg--Marquardt method (ILMM) for solving constrained nonsmooth equations. Since the orthogonal projection onto the feasible set may be computationally expensive, we propose a local ILMM  with feasible inexact projections. By using assumption of semi-smoothness and an error bound condition which is weaker than the standard full rank assumption, we establish its convergence results, including results on its rate. Finally, some computacional results are reported in order to illustrate the advantages of the new schemes.
\end{abstract}



\section{Introduction}\label{intro}
Here, we propose and investigue the inexact Levenberg--Marquardt method with feasible inexact projections for solving nonsmooth equations subject to a set convex, i.e., to solve the following problem: find $x\in \mathbb{R}^n$ such that
\begin{equation}\label{eq:prob}
\begin{aligned}
x\in C,\qquad f(x)=0,
\end{aligned}
\end{equation}
where  $C$ is a nonempty closed convex set contained in an open set $\Omega \subset \mathbb{R}^n$ and $f: \Omega \rightarrow \mathbb{R}^m$ is a locally Lipschitz continuous function.  Throughout this paper, we will assume that the solution set of \eqref{eq:prob}, denoted by $C^*$, is nonempty. The problem \eqref{eq:prob} has aroused the interest of many researchers since different problems can be written with \eqref{eq:prob}, for example, inequality feasibility problem and implicit complementarity problem, see \cite{PangQi1993,Pang1982}. Recently, in \cite{deOliveiraOFerreira2020} was proposed and analyzed a method for solving \eqref{eq:prob}. See also \cite{FacchiKanzon1997,KanzowPetra2004,KanPetra2007}. It is worth mentioning that, if $f$ is a continuously differentiable function, then \eqref{eq:prob} reduces to a constrained smooth equation, which has been addressed in several studies, and several methods have been proposed for solving it. See, for example, the exact/inexact Newton-like methods in \cite{mariniquasi2018,morini2016,GoncalvesGoncalves2Oliveira2021,GoncalvesOliveira2017,mariniquasi2018}, projected Levenberg--Marquardt-type methods in \cite{BehlingHaeserRamosSchonefeld2017,BehlingFischerHerrichIusemYe2014,BehlingRFischer2012}, and trust-region methods in \cite{BellaviMaria2004,Bellavia2012}.

Levenberg--Marquardt-type method has being a valuable tool for solving unconstrained/constrained equation, see \cite{Kenneth1944,Donald1963}. For solving a system of nonsmooth equations the exact version of this method is formulated as: given a current iterate $x_{k} \in \mathbb{R}^n$ and the parameter $\mu_{k}>0$, the next iterate is compute as $x_{k+1} = x_{k} + d_{k}$, where the vector $d_k$ is the solution of the liner system
\begin{equation}\label{uncon1}
(V_k^TV_k + \mu_kI_n)d = -V_k^Tf(x_k), \qquad V_k \in \partial f(x_k),
\end{equation}
where $I_n$ is the identity matrix in $\mathbb{R}^n\times \mathbb{R}^n$, $V_k := V_{x_k}$ is an element of the Clarke generalized Jacobian of $F$ at $x_k$. For the definition of the Clarke generalized Jacobian, see \cite{Clarke1990}. It is worth mentioning that, the solution of system \eqref{uncon1} is  equivalent the solution of the problem
\begin{equation}\label{uncon2}
\min_{d \in \mathbb{R}^n} \|f(x_k + V_kd)\|^2 +  \mu_k\|d\|^2,  \qquad V_k \in \partial f(x_k).
\end{equation}
The inexact versions of Levenberg--Marquardt-type method are more studied to the large-scale problems. In this case, at each iteration the system \eqref{uncon1} is not solved exactly, but only to within a certain tolerance, i.e, the vector $d_{k}$ is the solution of the system
\begin{equation}\label{unconstrained}
(V_k^TV_k + \mu_kI_n)d = -V_k^Tf(x_k) + r_k,
\end{equation}
where $r_k$ is a perturbation vector and measures how inexactly.

In this study, we propose a new scheme for solving \eqref{eq:prob}, which we refer to as ILMM-IP. Basically, the proposed method combines the inexact Levenberg--Marquardt method with a procedure to obtain a feasible inexact projection onto $C$, thus ensuring the viability of the iterates. The concept of feasible inexact projection was introduced in \cite{deOliveiraFerreiraSilva2018}, which also accepts an exact projection that can be adopted when it is easily obtained. For instance, exact projections onto a box constraint or Lorentz cone can be easily obtained; see \cite[p. 520]{NocedalWright2000} and \cite[Proposition 3.3]{FukushimaTseng2002}, respectively. It is noteworthy that a feasible inexact projection can be computed using any method that minimizes a quadratic function on a closed convex set efficiently by introducing a suitable error criterion. For instance, if the set $C$ is polyhedral, then some iterations of an interior point method or active set method can be performed to obtain a feasible inexact projection, see \cite{NicholasPhilippe2002,NocedalWright2000,Robert1996}. When $C$ is a simple convex compact set, a similar scheme is adopted \cite{GoncalvesOliveira2017,MaxJefferson2017,lan2016}, which uses the conditional gradient method to obtain a feasible inexact projection.

From the theoretical viewpoint, i.e., the local convergence of the proposed method as well as results on its rate are established by using assumption of semi-smoothness and an error bound condition. Specifically, let $\{x_k\}$ be the sequence generated by the method and $dist_{C^*}(x)$ the distance from $x$ to the solution set $C^*$. We show that the sequence $\{dist_{C^*}(x_k)\}$ converges to zero with a better rate than the linear. Moreover, we also deduce the convergence rate for the sequence $\{x_k\}$. To assess the practical behavior of the new method, some computacional results are reported. In particular, we compare the performance of the proposed method with that of the inexact Levenberg--Marquardt method using feasible exact projections.

The outline of this paper is as follows. In Section~\ref{sec:NotDef}, we present the notations and som technical definitions used herein. In Section~\ref{sec:condGmet3}, we describe the ILMM-IP and examine its local convergence analysis. Some preliminary numerical experiments for the proposed schemes are reported in Section~\ref{sec:CompResu}. Finally, some concluding remarks are given in Section~\ref{sec:Conclusions}.

\section{Notations and definitions}\label{sec:NotDef}
The inner product and its associated Euclidean norm in $\mathbb{R}^n$ are denoted by $\langle\cdot,\cdot\rangle$ and $\| \cdot \|$, respectively. The closed ball centered at $x$ with radius $\delta$ is denoted by $B_{\delta}(x) := \{ y \in \mathbb{R}^n \ : \ \| y - x\| \leq \delta\}$. We define the distance from $x$ to the solution set $C^*$ by
\begin{equation}\label{dist}
dist_{C^*}(x):= \inf_{y \in C^*} \| y - x \|,
\end{equation}
We also represent by $\bar{x}$ a point in $C^*$ which realizes such distance, i.e.
\begin{equation}\label{dist1}
\|x - \bar{x} \| := dist_{C^*}(x).
\end{equation}

\begin{definition}\label{def:NormOper}
The norm of a mapping $T: \mathbb{R}^n \to \mathbb{R}^m$ is define by
$$
\|T\|:= \sup \{\|Tx\|:~ x \in \mathbb{R}^n, \|x\| = 1\}.
$$
\end{definition}

In the following, we define the concept of a locally Lipschitz continuous function, which is crucial in our study.

\begin{definition}\label{def:FunLip}
A function $f:\Omega \subset \mathbb{R}^n \rightarrow \mathbb{R}^m$ regarded as locally Lipschitz continuous if for each $x\in\Omega$ there exist constants $L, \delta >0$ such that
$$
\|f(y) - f(z)\|\leq L \|y - z\| \quad \forall ~y, z \in B_\delta(x).
$$
\end{definition}
\begin{remark}
According to the Rademacher theorem, see \cite[Theorem~2, p.~81]{Evans1992}, locally Lipschitz continuous functions are differentiable almost everywhere.
\end{remark}

Next, we define the Clarke generalized Jacobian of a function, which has appeared in \cite{Clarke1990}. This Jacobian requires only the local Lipschitz continuity of the function $f$ and its well-definedness is ensured by the Rademacher theorem.

\begin{definition}\label{def:JacClarke}
The Clarke-generalized Jacobian of a locally Lipschitz continuous function $f$ at $x$ is a set-valued mapping $\partial f: \mathbb{R}^n \rightrightarrows \mathbb{R}^m$ defined as
$$
\partial f(x) := \mbox{co}\left\{H \in \mathbb{R}^{m\times n}:~ \exists \, \{x_k\} \subset \mathcal{D}_f, \lim_{k \to +\infty} x_k = x,\, H = \lim_{k \to + \infty}f'(x_k)\right\},
$$
where ``\mbox{co}'' represents the convex hull, $\mathbb{R}^{m\times n}$ the set comprising all $m\times n$ matrices, and $\mathcal{D}_f$ the set of points at which $f$ is differentiable.
\end{definition}

\begin{remark}
It is noteworthy that if $f$ is continuously differentiable at $x$, then $\partial f(x) = \{f'(x)\}$. Otherwise, $\partial f(x)$ may contain other elements different from $f'(x)$, even if $f$ is differentiable at $x$, see \cite[Example 2.2.3]{Clarke1990}. Furthermore, the Clarke-generalized Jacobian is a subset of $\mathbb{R}^{m\times n}$ that is typically nonempty, convex, and compact. In addition, the set-valued mapping $\partial f$ is closed and upper semi-continuous; see \cite[Proposition 2.6.2, p.~70]{Clarke1990}.
\end{remark}

We finish this section with another important result about the set-valued mapping $\partial f$.

\begin{proposition}\label{Prop:boundV}
The set-valued mapping $\partial f(\cdot): \mathbb{R}^n \rightrightarrows \mathbb{R}^m$ is locally bounded; that is, for all $\delta >0$, there exists a $L > 0$ such that for all $y \in B_{\delta}(x)$ and $V \in \partial f(y)$, $\|V\| \leq L$ holds.
\end{proposition}
\begin{proof}
For simplicity, we define the auxiliary set $\partial_B f(x)$ as follows
$$
\partial_B f(x) := \left\{H \in \mathbb{R}^{m\times n}:~ \exists \, \{x_k\} \subset \mathcal{D}_f, \lim_{k \to +\infty} x_k = x,\, H = \lim_{k \to + \infty}f'(x_k)\right\}.
$$
Since $V \in \partial f(y)$, there exist $H_1, \ldots, H_q \in  \partial_B  f(y)$   and $a_1, \ldots,  a_q\in [0,  1]$  such that $V = \sum_{\ell=1}^{q}a_\ell H_\ell$ and   $\sum_{\ell=1}^{m} a_\ell = 1$. On the other hand, because $H_1, \ldots, H_q \in  \partial_B  f(y)$, there exists $\{y_{k, \ell}\} \subset  B_\delta(y)\cap {\cal D}_f$ with $\lim_{k \rightarrow +\infty} y_{k, \ell} = y$ such that $V = \sum_{\ell=1}^{q}a_\ell\lim_{k\rightarrow +\infty}f'(y_{k, \ell})$. Since $\{y_{k, \ell}\} \subset  {\cal D}_f$, i.e., $f$ is differentiable at $y_{k, \ell}$ for each $\ell = 0, 1, \cdots, q$, we have  $f'(y_{k,\ell}, v) = f'(y_{k,\ell})v$, where $f'(y_{k,\ell}, v)$ is directional derivative of $f$ in the direction of vector $v$ at the point $y_{k,\ell}$. Using that $f'(y_{k,\ell}, v) = f'(y_{k,\ell})v$ and that $f$ is locally Lipschitz continuous, we obtain that
$$
\|f'(y_{k,\ell})v\| = \left\|\lim_{t \to 0^{-}} \dfrac{f(x + tv) - f(x)}{t}\right\| \leq L \|v\|.
$$
Hence, from Definition~\ref{def:NormOper}, we have $\|f'(y_{k,\ell})\| \leq L$. Now, using properties of the norm and the fact that $\sum_{\ell=1}^{m} a_\ell = 1$, we conclude
$$
\left\|V\right\| = \left\| \sum_{\ell=1}^{q} a_\ell \lim_{k\rightarrow +\infty}f'(y_{k, \ell}) \right\| \leq \sum_{\ell=1}^{q} \alpha_\ell  \lim_{k \rightarrow +\infty} \left\|f'(y_{k, \ell})  \right\| \leq L,
$$
which is the desired inequality.
\end{proof}
\section{Inexact Levenberg-Marquardt method with feasible inexact projections}\label{sec:condGmet3}

In this section, we propose and analyze a local inexact Levenberg-Marquardt method with feasible inexact projections (ILMM-IP) to solve problem \eqref{eq:prob}.
\begin{definition}\label{def:IP}
Let $x\in \mathbb{R}^n$ and $\epsilon \geq 0$ be given. We say that $ P_{C}(x,\epsilon)$ is an  $\epsilon$--projection of $x$ onto $C$ when
\begin{equation*}\label{eq:iproj}
P_{C}(x,\epsilon) \in C \quad \text{and} \quad \langle x - P_{C}(x,\epsilon), y - P_{C}(x,\epsilon) \rangle \leq \epsilon, \quad \forall y \in C.
\end{equation*}
\end{definition}

Next, we indicate some important information about the $\epsilon$--projection.

\begin{remark}\label{rem:projin}
If $\epsilon = 0$, then $P_C(x, 0)$ corresponds to the orthogonal projection of $x$ onto $C$, which will be denoted, simply, by $P_C(x)$. On the other hand, $P_C(x,\epsilon)$ is an $\epsilon$--projection of $x$ onto $C$ in the sense of Definition~\ref{def:IP} for any $\epsilon > 0$. In the case that orthogonal projection onto $C$ neither has a closed-form nor can be easily computed, an $\epsilon$--projection of $x$ onto $C$ can be obtained by means of an iterative method applied to solve the projection problem $\min_{y \in C}(0,5\|y - x\|^2)$. For example, if $C$ is bounded, one can use the Frank–Wolfe method \cite{fw1956}, to obtain an inexact projection in the sense of Definition~\ref{def:IP}. In particular, given $z_t \in C$, the t-th step of the Frank–Wolfe method first finds $w_t$ as a minimum of the linear function $\langle x - z_t , w_t - z_t\rangle \leq \epsilon$ over $C$ and then set $z_{t+1} = (1 - \alpha_t)z_t + \alpha_t w_t$ for some $\alpha_t \in [0,1]$.
\end{remark}

In what follows, we establish  a similar  property for the operator $P_{C}(\cdot,\cdot)$ whose proof can be found in \cite{GoncalvesOliveira2021}.

\begin{proposition}\label{prop:errproj}
For any  $x,y \in \mathbb{R}^n$ and $\epsilon\geq 0$, we have
\[
\|P_C(x,\epsilon) - P_C(y)\| \leq \|x-y\| + \sqrt{\epsilon}.
\]
\end{proposition}

We are now ready to formally described the inexact Levenberg-Marquardt method with feasible inexact projections.\\

\hrule
\begin{algorithm}  \label{Alg:NNM}
{\vspace{0.2cm}\bf ILMM-IP \vspace{0.3cm}}
\hrule
\vspace{0.1cm}
\begin{description}
\vspace{.5 cm}
\item[\bf Step 0.] Let $\eta \geq 1$, $\theta > 0$ and $\{\theta_k\}\subset[0,\theta)$ be given. Chose $x_0\in C$ and a parameter $\sigma \in (0,1)$. Set $k=0$.
\item[\bf Step 1.] If $f(x_k)=0$, then {\bf stop}.
\item[\bf Step 2.] Select an element $V_k \in \partial f(x_k)$ and set $\mu_k = \eta \|V_k^Tf(x_k)\|^\sigma$. Take a residual control $\zeta_k > 0$ and compute an approximate solution $d_k\in \mathbb{R}^n$  of the system
\begin{equation*}\label{A1:s1}
(V_k^T V_k + \mu_k I_n)d = -V_k^T f(x_k) + r_k,
\end{equation*}
such that
\begin{equation}\label{eq:10}
\|r_k\| \leq \zeta_k.
\end{equation}
\item[\bf Step 2.] Define $\epsilon_k := \theta_k^2 \| d_k \|^2$. Compute $P_{C}(x_k + d_k, \epsilon_k)$, an $\epsilon_k$--projection of $x_k + d_k$ onto $C$, and set
\begin{equation}\label{eq:lm3}
x_{k+1} := P_{C}(x_k + d_k, \epsilon_k).
\end{equation}
\item[\bf Step 3.] Set $k\gets k+1$, and go to \textbf{Step~1}.
\vspace{.5 cm}
\end{description}
\hrule
\end{algorithm}
\vspace{0.3cm}
\noindent
\begin{remark}
In ILMM-IP, we first verify whether the current iterate $x_k$ is a solution of \eqref{eq:prob}; otherwise, we select $V_k \in \partial f(x_k)$, $\mu_k > 0$ and $\zeta_k > 0$ to satisfy the criterion~\eqref{eq:10}. In our algorithm, we taking $\mu_{k} = \eta \|V_k^Tf(x_k)\|^\sigma$ for every $k\geq 0$ and $\sigma \in (0,1) $, which turns the problem \eqref{uncon1} into a strongly convex one and hence it possesses a unique solution. Suggestions of different regularization parameters have been discussed, see for example, \cite{Fan2013,Zhang2003}. As the point $x_k + d_k$ can be infeasible for the set of constraints $C$, ILMM uses a procedure to obtain a feasible inexact projection; consequently, the new iterate $x_{k+1}$ belongs to $C$. In particular, $x_{k+1}$ satisfies for $k\geq0$ the inequality
$$
\langle x - x_{k+1}, y - x_{k+1}\rangle \leq \epsilon, \qquad \forall~y \in C.
$$
See Remark~\ref{rem:projin} for some comments regarding our concept of feasible inexact projection and the method to compute it. Finally, the choose of tolerance $\theta_k$ is important for obtaining the local convergence of the ILMM-IP.
\end{remark}

In order to analyze the local convergence of the ILMM-IP, the following assumptions are made throughout this subsection.
\begin{itemize}
\item [\textbf{(H0)}] $C^* \neq \emptyset$ and let $x_* \in C^*$ be an arbitrary element of solution set.\\
\item [\textbf{(H1)}] There exist constants $\delta_1$, $\tau > 0$ and $0 \leq p \leq 1$ such that
$$
\|f(y) - f(x) - V(y - x)\| \leq \tau \|x - y\|^{1 + p},\quad \forall~ V \in \partial f(x), \quad \forall ~ x, y \in B_{\delta_1}(x^*).
$$
\item [\textbf{(H2)}] There exist $\omega, \delta_2 > 0$ such that $\|f(x)\|$ provides a local error bound on $B_{\delta_2}(x_*)$, i.e.,
$$
\omega\, dist_{C^*}(x) \leq \|f(x)\|, \qquad \forall~ x\in B_{\delta_2}(x^*).
$$
\end{itemize}

\begin{lemma}\label{Lemma1}
Suppose that assumptions (H0) to (H2) hold. Let $r_k$ the residual vector and consider $\delta := \min\{\delta_1, \delta_2, 2(\omega^2/L\tau)^{1/p}\}$. If $x_k \in B_{\delta/2}(x_*)$, then there exist constants $c_1, c_2 > 0$ such that
\begin{equation}\label{eq:8}
\|d_k\| \leq c_1 \,dist_{C^*}(x_k) + \dfrac{\|r_k\|}{\mu_k},
\end{equation}
\begin{equation}\label{equ:8}
\|V_kd_k  + f(x_k)\| \leq c_2 \,dist_{C^*}(x_k)^{1 + \frac{\sigma}{2}} + \|V_k\|\dfrac{\|r_k\|}{\mu_k}.
\end{equation}
\end{lemma}
\begin{proof}
Since $x_k \in B_{\delta/2}(x_*)$ from the triangular inequality, \eqref{dist} and \eqref{dist1}, we obtain that
$$
\|\bar{x}_k - x_*\| \leq \|\bar{x}_k - x_k\| + \|x_* - x_k\| \leq \|x_* - x_k\| + \|x_* - x_k\| \leq \delta,
$$
i.e., $\bar{x}_k \in B_{\delta}(x_*)$ and, consequently $\bar{x}_k \in B_{\delta_1}(x_*)$. On the other hand, because $f(\bar{x}_k) = 0$, we have
\begin{align*}
f(x_k)^TV_k(x_k - \bar{x}_k) &= f(x_k)^T[f(\bar{x}_k) + V_k(x_k - \bar{x}_k)] \nonumber\\
= f(x_k)^Tf(x_k) - f(x_k)^T[f(x_k) - f(\bar{x}_k) - V_k(x_k - \bar{x}_k)],
\end{align*}
where $V_k \in \partial f(x_k)$. Taking the norm on both sides of the last equality and using the properties of the norm, we have
\begin{equation}\label{equ:4}
\|f(x_k)^TV_k\|\|x_k - \bar{x}_k\| \geq \|f(x_k)\|^2 - \|f(x_k)\|\|f(x_k) - f(\bar{x}_k) - V_k(x_k - \bar{x}_k)\|.
\end{equation}
Since $x_k \in B_{\delta/2}(x_*)$ and $f(\bar{x}_k) = 0$, we can use assumption (H2), Definition~\ref{def:FunLip} and \eqref{dist1} to conclude that
\begin{equation}\label{equ:5}
\omega \|x_k - \bar{x}_k\| \leq \|f(x_k)\| \leq L\|x_k - \bar{x}_k\| \leq L\|x_k - x_*\|.
\end{equation}
Using assumption (H1), \eqref{equ:5} and the fact that $x_k \in B_{\delta/2}(x_*)$, we have \eqref{equ:4} reduces to
\begin{align*}
\|f(x_k)^TV_k\|\|x_k - \bar{x}_k\| &\geq \left[\omega^2 - L\tau \|x_k - x_*\|^p\right]\|x_k - \bar{x}_k\|^2 \geq \left[\omega^2 - L\tau \left(\dfrac{\delta}{2}\right)^p\right]\|x_k - \bar{x}_k\|^2.
\end{align*}
Because $\delta \leq 2(\omega^2/L\tau)^{1/p}$ the last inequality reduces to $\|f(x_k)^TV_k\| \geq \hat{c}\|x_k - \bar{x}_k\|$, where $\hat{c} := [\omega^2 - L\tau (\delta/2)^p]$. Now, note that
\begin{equation}\label{equ:6}
\mu_k = \eta \|V_k^Tf(x_k)\|^\sigma \geq \eta \hat{c}^\sigma \|x_k - \bar{x}_k\|^\sigma.
\end{equation}
On the other hand, using the properties of the norm, Proposition~\ref{Prop:boundV}, that $f(\bar{x}_k) = 0$ and Definition~\ref{def:FunLip}, we have
\begin{equation}\label{equ1}
\mu_k = \eta \|V_k^Tf(x_k)\|^\sigma \leq \eta L^\sigma \|f(x_k) - f(\bar{x}_k)\|^\sigma \leq \eta L^{2\sigma}\|x_k - \bar{x}_k\|^\sigma.
\end{equation}
Let $\bar{d}_k$ be solution of problem~\eqref{uncon2}, since $f(\bar{x}_k) = 0$, then by assumption (H1) and \eqref{equ:6}, we obtain that
\begin{align*}
\|\bar{d}_k\|^2 &\leq \dfrac{1}{\mu_k} \left[\|V_k \bar{d}_k + f(x_k)\|^2 + \mu_k\|\bar{d}_k\|^2\right] \\&\leq \dfrac{1}{\mu_k}\left[\|V_k(\bar{x}_k - x_k) + f(x_k)\|^2 + \mu_k\|\bar{x}_k - x_k\|^2\right]\nonumber \\
& \leq \dfrac{1}{\mu_k}\left[\|\tau^2 \|\bar{x}_k + x_k\|^{2 + 2p} + \mu_k\|\bar{x}_k - x_k\|^2\right] \\& \leq \left(\dfrac{\tau^2}{\eta \hat{c}^\sigma} \|\bar{x}_k - x_k\|^{2p - \sigma} + 1 \right)\|\bar{x}_k - x_k\|^2.
\end{align*}
Since $x_k \in B_{\delta/2}(x_*)$ and $\|\bar{x}_k - x_k\| \leq \|x_* - x_k\|$, we have
\begin{equation}\label{equ:7}
\|\bar{d}_k\| \leq \sqrt{\left(\dfrac{\tau^2}{\eta \hat{c}^\sigma} \left(\dfrac{\delta}{2}\right)^{2p - \sigma} + 1 \right)}\|\bar{x}_k - x_k\|.
\end{equation}
Moreover by \eqref{unconstrained}, follows that
\begin{align*}
d_k & = -(V_k^TV_k + \mu_k I_n)^{-1}V_k^T f(x_k) + (V_k^TV_k + \mu_kI_n)^{-1}r_k = \bar{d}_k + (V_k^TV_k + \mu_kI_n)^{-1}r_k.
\end{align*}
Taking the norm on both sides of the last equality, using the triangular inequality, \eqref{equ:7} and \eqref{dist1}, we conclude that
\begin{align*}
\|d_k\| &\leq \|\bar{d}_k\| + \|(V_k^TV_k + \mu_kI_n)^{-1}\|\|r_k\|\leq \sqrt{\left(\dfrac{\tau^2}{\eta \hat{c}^\sigma} \left(\dfrac{\delta}{2}\right)^{2p - \sigma} + 1 \right)}dist_{C^*}(x_k) + \dfrac{\|r_k\|}{\mu_k}.
\end{align*}
This implies that \eqref{eq:8} holds with $c_1 = \sqrt{(\tau^2/\eta \hat{c}^\sigma)(\delta/2)^{2p - \sigma} + 1}$. We proceed to prove the inequality in \eqref{equ:8}. Considering the left-hand of \eqref{equ:8}, we have
\begin{align}\label{equ4}
\|V_k d_k + f(x_k)\| &= \|V_k \left[\bar{d}_k + (V_k^TV_k + \mu_k I_n)^{-1}r_k \right] + f(x_k)\|\nonumber\\
& \leq \|V_k \bar{d}_k + f(x_k)\| + \|V_k\|\|(V_k^TV_k + \mu_k I_n)^{-1}\|\|r_k\|\nonumber\\
& \leq \|V_k \bar{d}_k + f(x_k)\| + \|V_k\|\dfrac{\|r_k\|}{\mu_k}.
\end{align}
Since $\bar{d}_k$ is solution of problema~\eqref{uncon2} and $f(\bar{x}_k) = 0$, using assumption (H1), and \eqref{equ1}, we obtain that
\begin{align}\label{equ:9}
\|V_k\bar{d}_k + f(x_k)\|^2 & \leq \|V_k\bar{d}_k + f(x_k)\|^2 + \mu_k \|\bar{d}_k\|\nonumber\\ &\leq \|V_k(\bar{x}_k - x_k) + f(x_k)\|^2 + \mu_k \|\bar{x}_k - x_k\|^2 \nonumber\\
& \leq \tau^2 \|\bar{x}_k - x_k\|^{2 + 2p} + \eta L^{2\sigma}\|\bar{x}_k - x_k\|^{2 + \sigma}\nonumber\\ &\leq \left[\tau^2\left(\dfrac{\delta}{2}\right)^{2p-\sigma} +\eta L^{2\sigma}\right]\|\bar{x}_k - x_k\|^{2+\sigma}.
\end{align}
Therefore, extracting the square root on both sides of \eqref{equ:9} and combining with \eqref{equ4}, we obtain the desired result.
\end{proof}

In the following, we give an assumption about the residual vector $r_k$.

\begin{itemize}
\item [\textbf{(H3)}] Let $\sigma \in (0,1)$ and $\nu_k \subset \mathbb{R}_+$ for all $k = 0, 1, \ldots$. The residual vector $r_k$ satisfies
$$
\dfrac{\|r_k\|}{\mu_k}\leq \nu_k \, dist_{C^*}(x_k).
$$
\end{itemize}
In this paper, we choose $\nu_k \leq dist_{C^*}(x_k)^{\frac{\sigma}{2}}$ for all $k = 0, 1, \ldots$. If assumption (H3) holds, there exists $\delta_3 >0$ such that
\begin{equation}\label{eq:6}
dist_{C^*}(x_k) \leq \delta_3 \quad \Rightarrow \quad \dfrac{\|r_k\|}{\mu_k} \leq dist_{C^*}(x_k)^{1+ \frac{\sigma}{2}} \leq dist_{C^*}(x_k).
\end{equation}

\begin{lemma}\label{Lemm2}
Suppose that assumptions (H0) to (H3) hold and let $\{x_k\}$ be a sequence generated by Algorithm~\ref{Alg:NNM}. Let $\delta \leq \min\{\delta_1, \delta_2, \delta_3, 2(\omega^2/L\tau)^{1/p}\}$, where $\delta_1$, $\delta_2$ are defined in Lemma~\ref{Lemma1} and $\delta_3$ in assumption (H3). Then for any $x_k$, $x_k + d_k \in B_{\delta/2}(x_*)$, there exists $c_3 > 0$ such that
$$
dist_{C^*}(x_{k+1}) < c_3 \, dist_{C^*}(x_k)^{1+\frac{\sigma}{2}},
$$
for all $k = 0, 1, \ldots$.
\end{lemma}

\begin{proof}
Follows on from \eqref{dist} and \eqref{eq:lm3} that
\begin{equation}\label{eq:12}
dist_{C^*}(x_{k+1}) = dist_{C^*}\big(P_C(x_k + d_k, \epsilon_k)\big) = \inf_{x\in C^*}\|P_C(x_k + d_k, \epsilon_k) - x\|.
\end{equation}
Since $P_C(x) = x$ for each $x\in C$, we obtain from Proposition~\ref{prop:errproj} and \eqref{dist} that
\begin{align}\label{eq:2}
\inf_{x \in C^*}\|P_C(x_k + d_k, \epsilon_k) - P_C(x)\| &\leq \inf_{x \in C^*}(\|x_k + d_k - x\| + \sqrt{\epsilon_k})\nonumber\\
& = \sqrt{\epsilon_k} + \inf_{x \in C^*}\|x_k + d_k - x\|\nonumber\\
& = \sqrt{\epsilon_k} + dist_{C^*}(x_k + d_k).
\end{align}
Combining \eqref{eq:12} and \eqref{eq:2}, and using assumption (H2) and the fact that $\epsilon_k = \theta_k^2\|d_k\|^2$, we have
\begin{equation}\label{eq:3}
dist_{C^*}(x_{k+1}) \leq \sqrt{\epsilon_k} + dist_{C^*}(x_k + d_k) \leq \theta_k\|d_k\| + \dfrac{1}{\omega}\|f(x_k + d_k)\|.
\end{equation}
On the other hand by assumption (H1), we obtain that
\begin{equation}\label{eq:4}
\|f(x_k + d_k)\| - \|f(x_k) + V_kd_k\|\leq \|f(x_k + d_k) - f(x_k) - V_kd_k\| \leq \tau \|d_k\|^{1 + p},
\end{equation}
where $V_k \in \partial f(x_k)$. Combining \eqref{eq:3} and \eqref{eq:4}, we have
$$
dist_{C^*}(x_{k+1}) \leq \theta_k\|d_k\| + \dfrac{1}{\omega}(\tau \|d_k\|^{1 + p} + \|f(x_k) + V_kd_k\|).
$$
Now, we can use Proposition~\ref{Prop:boundV}, Lemma~\ref{Lemma1} and \eqref{eq:6} to conclude that
\begin{multline*}
dist_{C^*}(x_{k+1})  \leq \left[\theta_k (c_1 + 1)dist_{C^*}(x_k)^{- \frac{\sigma}{2}} + \dfrac{\tau (c_1 + 1)^{1 + p}}{\omega}dist_{C^*}(x_k)^{p - \frac{\sigma}{2}}\right.\\ \left.+\dfrac{c_2 + L}{\omega}\right]dist_{C^*}(x_k)^{1+\frac{\sigma}{2}}.
\end{multline*}
Finally, since $x_k \in B_{\delta/2}(x_*)$ and $\theta_k < \theta$. Then, the last inequality reduces to
$$
dist_{C^*}(x_{k+1}) < \left\{\left[\theta(c_1 + 1)\left(\dfrac{2}{\delta}\right)^p + \dfrac{\tau (c_1+1)^{1+p}}{\omega}\right]{\left(\dfrac{\delta}{2}\right)}^{p - \frac{\sigma}{2}}+ \dfrac{c_2 + L}{\omega} \right\} dist_{C^*}(x_k)^{1+\frac{\sigma}{2}},
$$
for all $k = 0, 1, \ldots$. Therefore, the proof of the lemma is complete with constant $c_3 = [\theta(c_1 + 1)(2/\delta)^p + \tau (c_1+1)^{1+p}/\omega]{(\delta/2)}^{p - \frac{\sigma}{2}} + (c_2 + L)/\omega $.
\end{proof}

In the next result, we show that $x_k, x_k + d_k \in B_{\delta/2}(x_*)$ if the starting point $x_0$ in the proposed algorithm is chosen sufficiently close to the solution set $C^*$.

\begin{lemma}\label{Lemma3}
Let $\hat{\delta} := \min\{1/2c_3^{-1/p}, \delta/ 2(c_1 + 2)[1 + (1+\theta)(c_1 + 1)\varsigma]\}$, where $\varsigma \geq \sum_{\ell = 0}^{\infty} (1/2)^{(1 + p)^{\ell} - 1}$. If $x_0 \in B_{\hat{\delta}}(x^*)\cap C$, then $x_k, x_k + d_k \in B_{\delta/2}(x_*)$ for every $k \geq 0$.
\end{lemma}
\begin{proof}
We will proceed by induction on $k$. We start with $k = 0$. By assumption, we have $x_0 \in B_{\hat{\delta}}(x^*)$. Since $\hat{\delta} \leq \delta/2$, we conclude that $x_0 \in B_{\delta/2}(x^*)$. Moreover, using the triangular inequality, \eqref{eq:8}, \eqref{eq:6} and \eqref{dist}, we obtain
$$
\|x_0 + d_0 - x_*\| \leq \|x_0 - x_*\| + \|d_0\| \leq \hat{\delta} + (c_1 + 1)\,dist_{C^*}(x_0) \leq \hat{\delta} + (c_1 + 1)\|x_0 - x_*\|
\leq (c_1 + 2)\hat{\delta}.
$$
Since, in particular, $\hat{\delta} \leq \delta/2(c_1 + 2)$, we conclude that $x_0 + d_0 \in B_{\delta/2}(x_*)$. Let $k \geq 0$ be arbitrarily given and assume that $x_{\ell}, x_{\ell} + d_{\ell} \in B_{\delta/2}(x_*)$ for all $\ell = 0, 1, \ldots, k$. Now, we procede to prove that $x_{k+1}, x_{k+1} + d_k \in B_{\delta/2}(x_*)$. By induction assumption $x_{\ell}, x_{\ell} + d_{\ell} \in B_{\delta/2}(x_*)$ for all $\ell = 0, 1, \ldots, k$, thus from Lemma~\ref{Lemm2} and \eqref{dist}, we have
\begin{align}\label{equ:1}
dist_{C^*}(x_{\ell}) & < c_3\cdot c_3^{1+p} \cdot c_3^{(1+p)^2}\cdot\ldots\cdot c_3^{(1+p)^{\ell - 1}} \,dist_{C^*}(x_0)^{(1 + p)^{\ell}}\nonumber \\
& \leq c_3^{\frac{1}{p}[(1+p)^{\ell} - 1]}\|x_0 - x_*\|^{(1 + p)^{\ell}} \leq c_3^{\frac{1}{p}[(1+p)^{\ell} - 1]}\hat{\delta}^{(1 + p)^{\ell}}.
\end{align}
On the other hand, using \eqref{eq:lm3} and Proposition~\ref{prop:errproj}, we find that
$$
\|x_{k + 1} - x_*\| = \|P_{C}(x_k + d_k, \epsilon_k) - P_C(x_*)\| \leq \|x_k + d_k - x_*\| + \sqrt{\epsilon_k}.
$$
By the triangle inequality and the facts that $\epsilon_k = \theta_k^2\|d_k\|^2$ and $\theta_k < \theta$, we obtain that
$$
\|x_{k+1} - x_*\| < \|x_k - x_*\| + (1 + \theta)\|d_k\|.
$$
Since $x_{\ell} \in B_{\delta/2}(x_*)$, for all $\ell = 0,1, \ldots, k$ and $x_0 \in B_{\hat{\delta}}(x_*)$, we can use the last inequality recursively together with \eqref{eq:8} and \eqref{eq:6} to conclude that
\begin{align}\label{equ:2}
\|x_{k+1} - x_*\| & < \|x_0 - x_*\| + (1 + \theta)\sum_{\ell = 0}^{k}\|d_{\ell}\| \nonumber \\ &\leq \hat{\delta} + (1 + \theta)(c_1 + 1)\sum_{\ell = 0}^k dist_{C^*}(x_{\ell})\nonumber \\
& \leq \hat{\delta} + (1 + \theta)(c_1 + 1)\sum_{\ell = 0}^{\infty} dist_{C^*}(x_{\ell}).
\end{align}
Now combining \eqref{equ:1} and \eqref{equ:2}, we obtain that
\begin{align}\label{equ3}
\|x_{k+1} - x_*\| & < \hat{\delta} + (1 + \theta)(c_1 + 1)\sum_{\ell = 0}^{\infty} c_3^{\frac{1}{p}[(1+p)^{\ell} - 1]}\hat{\delta}^{(1 + p)^{\ell}} \nonumber
\\ &\leq \hat{\delta} + (1 + \theta)(c_1 + 1)\hat{\delta}\sum_{\ell = 0}^{\infty} \left(\dfrac{1}{2}\right)^{(1 + p)^{\ell} - 1}\nonumber \\
& \leq  [1 + (1 + \theta)(c_1 + 1)\varsigma]\hat{\delta},
\end{align}
where the second inequality follows from definition $\hat{\delta}$. Therefore, $x_{k+1} \in B_{\delta/2}(x_*)$. It remains to prove that $x_{k+1} + d_{k+1} \in B_{\delta/2}(x_*)$. Since $x_{k+1} \in B_{\delta/2}(x_*)$, it follows from the triangle inequality, \eqref{eq:8} and \eqref{eq:6} that
$$
\|x_{k+1} + d_{k+1} - x_*\|\leq \|x_{k+1} - x_*\| + \|d_{k+1}\| \leq (c_1 + 2)\|x_{k+1} - x_*\|,
$$
which, combined with \eqref{equ3} and the definition of $\hat{\delta}$, yields
$$
\|x_{k+1} + d_{k+1} - x_*\| <  (c_1 + 2)[1 + (1 + \theta)(c_1 + 1)\varsigma]\hat{\delta} \leq \dfrac{\delta}{2}.
$$
This implies $x_{x+1} + d_{k+1} \in B_{\delta/2}(x_*)$ and then the proof is complete.
\end{proof}

We are now ready to prove the convergence of the sequences $\{dist_{C^*}(x_k)\}$ and $\{x_k\}$.
\begin{theorem}
Suppose that assumptions (H0) to (H3) hold and let $\{x_k\}$ be a sequence generated by Algorithm~\ref{Alg:NNM} with $x_0 \in B_{\hat{\delta}}(x_*)\cap C$. Let $\delta$ and $\hat{\delta}$ be the constants given in Lemmas~\ref{Lemma1} and \ref{Lemma3}, respectively. Then the sequence $\{dist_{C^*}(x_k)\}$ converges to zero at $(1 + \sigma/2)$ - order rate and the sequence $\{x_k\}$ converges to a point belonging to $C^*$.
\end{theorem}
\begin{proof}
The first part follows immediately from Lemmas~\ref{Lemm2} and \ref{Lemma3}. Now, we proceed to prove the second part. Since $\{dist_{C^*}(x_k)\}$ converges to zero and from Lemma~\ref{Lemma3}, $\{x_k\} \subset B_{\hat{\delta}}(x^*)\cap C \subset  B_{\delta/2}(x_*)\cap C$, it suffices to show that $\{x_k\}$ converges. Let us prove that $\{x_k\}$ is a Cauchy sequence. For this end, take $p, q \in \mathbb{N}$ with $p \geq q$. It follows from Proposition~\ref{prop:errproj}, triangular inequality, and the facts that $\epsilon_k = \theta_k^2 \|d_k\|^2$ and $\{x_k\} \subset C$ that
\begin{align*}
\|x_p - x_q\|& = \|P_C(x_{p-1} + d_{p-1}, \epsilon_{p-1}) - P_C(x_q)\| \nonumber\\
& \leq \|x_{p - 1} + d_{p-1} - x_q\| + \theta_{p-1}\|d_{p-1}\| \nonumber\\
& \leq \|x_{p - 1} - x_q\| + (1 + \theta_{p-1})\|d_{p-1}\|.
\end{align*}
Repeating the process above, we get
$$
\|x_p - x_q\| \leq (1 + \theta_q)\|d_q\| + \cdots + (1 + \theta_{p - 2})\|d_{p - 2}\| + (1 + \theta_{p-1})\|d_{p-1}\|,
$$
which, combined with the fact that $\theta_k < \theta$, for every $k \geq 0$, yields
$$
\|x_p - x_q\| < (1 + \theta)\sum_{\ell = q}^{p-1}\|d_{\ell}\|  \leq (1 + \theta)\sum_{\ell = q}^{\infty}\|d_{\ell}\|.
$$
On the other hand by \eqref{eq:8}, \eqref{eq:6}, \eqref{equ:1} and definition of $\hat{\delta}$, we have
$$
\|d_k\| \leq (c_1 + 1)\,dist_{C^*}(x_{\ell}) \leq (c_1 + 1)c_3^{\frac{1}{p}[(1+p)^{\ell} - 1]}\hat{\delta}^{(1+p)^{\ell}} \leq (c_1 + 1)\hat{\delta} \left(\dfrac{1}{2}\right)^{(1+p)^{\ell} - 1}.
$$
Combining the last two inequalities, we obtain that
\begin{align*}
\|x_p - x_q\| & < (1 + \theta)(c_1 + 1)\hat{\delta}\sum_{\ell = q}^{\infty} \left(\dfrac{1}{2}\right)^{(1+p)^{\ell} - 1} \\ &= (1 + \theta)(c_1 + 1)\hat{\delta}\left[\sum_{\ell = 0}^{\infty} \left(\dfrac{1}{2}\right)^{(1+p)^{\ell} - 1} - \sum_{\ell = 0}^{q-1} \left(\dfrac{1}{2}\right)^{(1+p)^{\ell} - 1}\right].
\end{align*}
Taking the limit in the last inequality as $q$ goes to $\infty$, we have $\|x_p - x_q\|$ goes to $0$. Therefore, $\{x_k\}$ is a Cauchy sequence and hence it converges. Let us say $\bar{x} = \lim_{k \to \infty}x_k$. Since $x_k \in C$, $\forall\, k$, and $C$ is closed, them $\bar{x} \in C$. Moreover, because $\omega \,dist_{C^*}(x_k) \leq \|f(x_k)\| \leq L \,dist_{C^*}(x_k)$, $f$ is continuous and $\{dist_{C^*}(x_k)\}$ converges to zero, we conclude that $\bar{x} \in C^*$.
\end{proof}

The following theorem proves the local convergence rate of the sequence $\{x_k\}$ generated by the ILMM-IP.

\begin{theorem}
Suppose the assumptions (H0) to (H3) hold and let $\{x_k\}$ be a sequence generated by Algorithm~\ref{Alg:NNM} with $x_0 \in B_{\hat{\delta}}(x_*)\cap C$. Let $\delta$ and $\hat{\delta}$ be the constants given in Lemmas~\ref{Lemma1} and \ref{Lemma3}. Then the sequence $\{x_k\}$ converges to $\bar{x}$ at $(1 + \sigma/2)$ - order rate.
\end{theorem}
\begin{proof}
It follows from Proposition~\ref{prop:errproj} that
$$
\|d_k\| = \|x_k + d_k - x_k\| \geq \|P_C(x_k + d_k, \epsilon_k) - P_C(x_k)\| - \sqrt{\epsilon_k}.
$$
Using \eqref{eq:lm3} and the facts that $\epsilon_k = \theta_k^2 \|d_k\|^2$ and $\theta_k < \theta$, we conclude that
$$
\|d_k\| \geq \|x_{k+1} - x_k\| - \theta \|d_k\|,
$$
Now let $\bar{x}_{k+1}$ be satisfying $dist_{C^*}(x_{k+1}) = \|x_{k+1} - \bar{x}_{k+1}\|$. Thus, from the previous inequality and \eqref{dist1}, we have
$$
(1 + \theta)\|d_k\| \geq \|x_{k} - \bar{x}_{k+1}\| - \|\bar{x}_{k+1} - x_{k+1}\| = dist_{C^*}(x_k) - dist_{C^*}(x_{k+1}).
$$
From Lemma~\ref{Lemm2}, we can get that $dist_{C^*}(x_{k+1}) < dist_{C^*}(x_k)/2$, when $k$ is sufficiently large. Hence, $\|d_k\| > dist_{C^*}(x_k) / 2(1 + \theta) $ and from \eqref{eq:8}, \eqref{eq:6} and Lemma~\ref{Lemm2}, we conclude that
\begin{equation}\label{equ:dis}
\|d_{k+1}\|  \leq (c_1 + 1)\,dist_{C^*}(x_{k+1}) < c_3(c_1 + 1)\,dist_{C^*}(x_{k})^{1 + \frac{\sigma}{2}} < c_3(c_1 + 1)2^{1 + \frac{\sigma}{2}}(1 + \theta)^{1 + \frac{\sigma}{2}}\|d_k\|^{1 + \frac{\sigma}{2}}.
\end{equation}
For sufficiently large $k$, we assume without loss of generality that the condition $c_3(c_1 + 1)2^{1 + \frac{\sigma}{2}}(1 + \theta)^{1 + \frac{\sigma}{2}}\|d_k\|^{\frac{\sigma}{2}} \leq 1/2$ holds. Therefore, we can get that $\|d_{k+1}\| < 1/2 \|d_k\|$ and consequently, we have
\begin{equation}\label{eq:7}
\|d_{k+j}\| < \dfrac{1}{2}\|d_{k+j-1}\| < \cdots < \left(\dfrac{1}{2}\right)^j\|d_k\|, \qquad j = 0, 1, \cdots.
\end{equation}
Using Proposition~\ref{prop:errproj} and the fact that $\epsilon_k = \theta_k^2\|d_k\|^2$, we obtain that
\begin{align*}
\|x_k - x_{x+l}\| & = \|P_C(x_k) - P_C(x_{k+l-1} + d_{k+l-1}, \epsilon_{k+l-1})\| \\
                  & \leq \|x_k - x_{k+l-1} + d_{k+l-1}\|  + \sqrt{\epsilon_{k+l-1}} \\
                  & \leq \|x_k - x_{k+l-1}\| + (1 + \theta_{k+l-1})\|d_{k+l-1}\|.
\end{align*}
Repeating the process above, we get
$$
\|x_k - x_{x+l}\| \leq (1 + \theta_k)\|d_k\| + \cdots + (1 + \theta_{k+l-2})\|d_{k+l-2}\| + (1 + \theta_{k+l-1})\|d_{k+l-1}\|,
$$
which, combined with the fact $\theta_k < \theta$, for every $k \geq 0$, and \eqref{eq:7}, yields
$$
\|x_k - x_{x+l}\| < (1 + \theta)\sum_{j = 0}^{l-1}\|d_{k+j}\| < (1 + \theta)\|d_k\|\sum_{j = 0}^{l-1}\left(\dfrac{1}{2}\right)^{j}.
$$
Since $\bar{x}$ is the limit point of $\{x_k\}$ taking the limit in the last inequality as $l$ goes to $\infty$, we obtain that
\begin{equation}\label{eq:9}
\|x_k - \bar{x}\| = \lim_{l \to \infty}\|x_k - x_{x+l}\|  \leq (1 + \theta)\|d_k\|\sum_{j = 0}^{\infty}\left(\dfrac{1}{2}\right)^{j}.
\end{equation}
Considering that $\sum_{j = 0}^{\infty} (1/2)^j = 2$, we conclude from \eqref{eq:9} that $\|x_k - \bar{x}\| \leq 2 (1 + \theta)\|d_k\|$. Therefore, from \eqref{eq:8}, \eqref{eq:6}, \eqref{equ:dis} and \eqref{dist1}, we have
\begin{equation*}
\|x_{k+1} - \bar{x}\| \leq 2 (1 + \theta)\|d_{k+1}\| < c_3(c_1 + 1)[2(1+ \theta)]^{2 + \frac{\sigma}{2}}\|d_k\|^{1 + \frac{\sigma}{2}} \leq c_3[2(c_1 + 1)(1+ \theta)]^{2 + \frac{\sigma}{2}}\|x_k - \bar{x}\|^{1 + \frac{\sigma}{2}},
\end{equation*}
which implies that $\{x_k\}$ converges to $\bar{x}$ at $(1+\sigma/2)$ - order rate.
\end{proof}

\section{Computacional results}\label{sec:CompResu}

Here, we present some computacional results to assess the practical behavior of inexact Levenberg-Marquardt method with exact
projections (ILMM-EP) and inexact Levenberg-Marquardt method with feasible inexact projections (ILMM-IP). Specifically, we consider one classe of constrained nonsmooth equations, i.e,  of the form \eqref{eq:prob}. We worked with one class of medium- and large-scale problems called CAVEs (constrained absolute value equations). It is worth mentioning that in \cite{DeOliveiraFerreira2020} a version of inexact Newton method was used for solving CAVEs. The Constrained absolute value equation (CAVE) is described as
$$
\mbox{find} \quad x \in C \quad \mbox{such that}\quad Ax - |x| = b,
$$
where $C := \{x \in \mathbb{R}^n:~ \sum_{i = 1}^{n}x_i \leq d, \, x_i \geq 0, \, i = 1,2, \ldots,n\}$, $A \in \mathbb{R}^{n\times n}$, $b \in \mathbb{R}^n \equiv \mathbb{R}^{n\times 1}$, and $|x|$ denotes the vectors whose $i$-th component is equal to $|x_i|$.

In our implementation, the CAVEs have been generated randomly. We used the Matlab routine \textit{sprand} to construct matrix $A$. In particular, this routine generates a sparse matrix with predefined dimension, density, and singular values. Initially, we defined the dimension $n$ and randomly generated the vector of singular values from a uniform distribution on $(0, 1)$. To ensure that $\|A^{-1}\|< 1/3$, i.e., so that the assumptions of \cite[Theorem 2]{BelloCruz2016} are fulfilled, we rescale the vector of singular values by multiplying it by $3$ divided by the minimum singular value multiplied by a random number in the interval $(0, 1)$. To generate the vector $b$ and the constant $d$, we chose a random solution $x_*$ from a uniform distribution on $(0.1, 100)$ and computed $b = Ax_* - |x_*|$ and $d = \sum_{i = 1}^{n}(x_*)_i$, where $(x_*)_i$ denotes the $i$-th component of the vector $x_*$. In both methods, $x_0 = (d/2n, d/2n,\ldots, d/2n)$ was defined as the starting point, the initialization data $\theta$ was taken equal to $10^{-2}$. We stopped the execution of Algorithm~\ref{Alg:NNM} at $x_k$, declaring convergence if $\|Ax_k - |x_k| - b\| < 10^{-6}$. In case this stopping criterion was not respected, the method stopped when a maximum of $100$ iterations had been performed.  The procedure to obtain feasible projections used in our implementation was the \textit{CondG Procedure}; see, for example, \cite{OliveiraFerreiraSilva2019, GoncalvesOliveira2017}. In particular, this procedure stopped when either the stopping criterion, i.e., the condition $\langle x_k - x_{k+1}, y - x_{k+1}\rangle \leq  \epsilon$ was satisfied for all $y \in C$ and $k = 0, 1, \ldots$ or a maximum of $100$ iterations was performed. For this class of problems, an element of the Clarke generalized Jacobian at $x$ (see \cite{BelloCruz2016,Mangasarian2009}) is given by
$$
V = A -  \mbox{diag}(\mbox{sgn}(x)), \qquad x \in  \mathbb{R}^n,
$$
where $\mbox{diag}(\alpha_i)$ denotes a diagonal matrix with diagonal elements $\alpha_1, \alpha_2, \ldots, \alpha_n$ and $\mbox{sgn}(x)$ denotes a vector with components equal to $-1$, $0$, or $1$ depending on whether the corresponding component of the vector $x$ is negative, zero, or positive.
The ILMM-EP and ILMM-IP requires that the linear system $(V_k^TV_k + \mu_k I_n)d + V_k^Tf(x_k) - r_k = 0$ to be solved approximately, in the sense of \eqref{eq:10} be satisfied. Matlab has several iterative methods for solving linear equations. For our class of problems, the routine \textit{bicgstab} was the most efficient; thus, in all tests, we used \textit{bicgstab - BiConjugate Gradients Stabilized Method} as an iterative method to solve linear equations approximately. The numerical results were obtained using Matlab version R2016a on a 2.5~GHz Intel\textregistered\ Core\texttrademark\ i5 2450M computer with 6~GB of RAM and Windows 7 ultimate system.

Table~\ref{tab:res} displays the numerical results obtained for the test set proposed. The methods were compared on the total number of iterations (It) and CPU time in seconds (Time).

Observe that as for the CPU time the ILMM-IP behaved better than the ILMM-EP for this problem set. In particular, this becomes more evident for problems of large-scale. Regarding the total number of iterations, both methods behaved quite similarly.

In summary, the numerical experiments indicate that the ILMM-IP seems to be reliable and competitive for solving medium- and large-scale constrained nonsmooth mainly when the orthogonal projection onto the feasible set cannot be easily computed.

\begin{table}
\centering
\caption{Comparison of ILMM-EP and ILMM-IP for solving CAVEs.}\label{tab:res}
\begin{tabular}{rrrcccc}
\hline
\ & \ & \multicolumn{2}{c}{ILMM-EP} & \multicolumn{2}{c}{ILMM-IP} \\
$n$ & $m$ &   It & Time & It & Time \\ \hline
100 & 100 &  6 & 0.15 & 7 & 0.13 \\ \hline
500 & 500 &  6 & 0.56 & 6 & 0.49  \\ \hline
1000 & 1000 &7 & 3.38 & 7 & 2.65  \\ \hline
5000 & 5000 &8 & 162.90& 8 & 152.13  \\ \hline
\end{tabular}
\end{table}

\section{Conclusions}\label{sec:Conclusions}
This paper proposed and analyzed a ILMM-IP. It basically, consists of combining  the inexact Levenberg-Marquart method step with feasible inexact projections. The local convergence as well as results on its rate were established under assumption of semi-smoothness and an error bound condition, which is weaker than the standard full-rank condition of Clarke-generalized Jacobian of $f$ at $x$. Finally, some numerical experiments were carried out in order to illustrate the numerical behavior of proposed method. In particular, they indicate that the ILMM-IP represent a useful tool for solving medium- and large-scale constrained nonsmooth equations mainly when the orthogonal projection onto the feasible set cannot be easily computed.

\end{document}